\documentclass[12pt]{amsart}
\usepackage{amsfonts}
\usepackage{amsmath,amssymb,amsthm}
\usepackage{color}
\input epsf.sty

\newcommand{\ds}{\displaystyle}
\newtheorem{thm}{Theorem}[section]
\newtheorem{cor}[thm]{Corollary}

\newtheorem{prop}[thm]{Proposition}
\theoremstyle{definition}

\theoremstyle{remark}
\newtheorem{rem}[thm]{Remark}

\begin{document}
\title[POINTS, WHOSE PEDAL TRIANGLES ARE SIMILAR TO THE GIVEN TRIANGLE]
{POINTS, WHOSE PEDAL TRIANGLES ARE SIMILAR TO THE GIVEN TRIANGLE}
\author{Georgi Ganchev$^1$, \; Gyulbeyaz Ahmed$^2$ \; and \; Marinella Petkova$^3$}%

\address{1 Bulgarian Academy of Sciences, Institute of Mathematics and Informatics,
Acad. G. Bonchev Str. bl. 8, 1113 Sofia, Bulgaria}%
\email{1 ganchev@math.bas.bg}%
\address{2,3 High school ''P. R. Slaveykov'' Kardjali, Bulgaria}%

\subjclass[2000]{Primary 53A05, Secondary 53A10}%
\keywords{Brocard points, pedal triangles}%

\begin{abstract}
We study the eleven points in the plane of a given triangle, whose
pedal triangles are similar to the given one. We prove that the six
points whose pedal triangles are positively oriented, lie on a single
circle, while the five points, whose pedal triangles are negatively
oriented, lie on a common straight line.
\end{abstract}

\maketitle
\section{Preliminaries}

First, following mainly \cite{G} we give some notions and facts,
which we use further in the paper (see also \cite{GN}).

\subsection{The map $f$}
Let $ABC$ be a triangle with sides $BC=a$, $CA=b$, $AB=c$, and angles
$\angle A=\alpha, \angle B = \beta, \angle C = \gamma$, which we call the
\emph{basic triangle}. For an arbitrary point $M$ in the plane of
$\triangle ABC$ we denote the segments $MA=x$, $MB=y$ and $MC=z$.

Any circle inversion $\varphi(M,r)$ with center $M$ and an arbitrary radius
$r$ maps the points $ABC$ into the points $A_1B_1C_1$ (Figure 1),
respectively, so that
$$B_1C_1=\frac{r^2}{xyz}\,ax, \quad C_1A_1=\frac{r^2}{xyz}\,by, \quad
A_1B_1=\frac{r^2}{xyz}\,cz.\leqno(1)$$

\begin{figure}[h]\center\epsfysize=5cm\epsffile{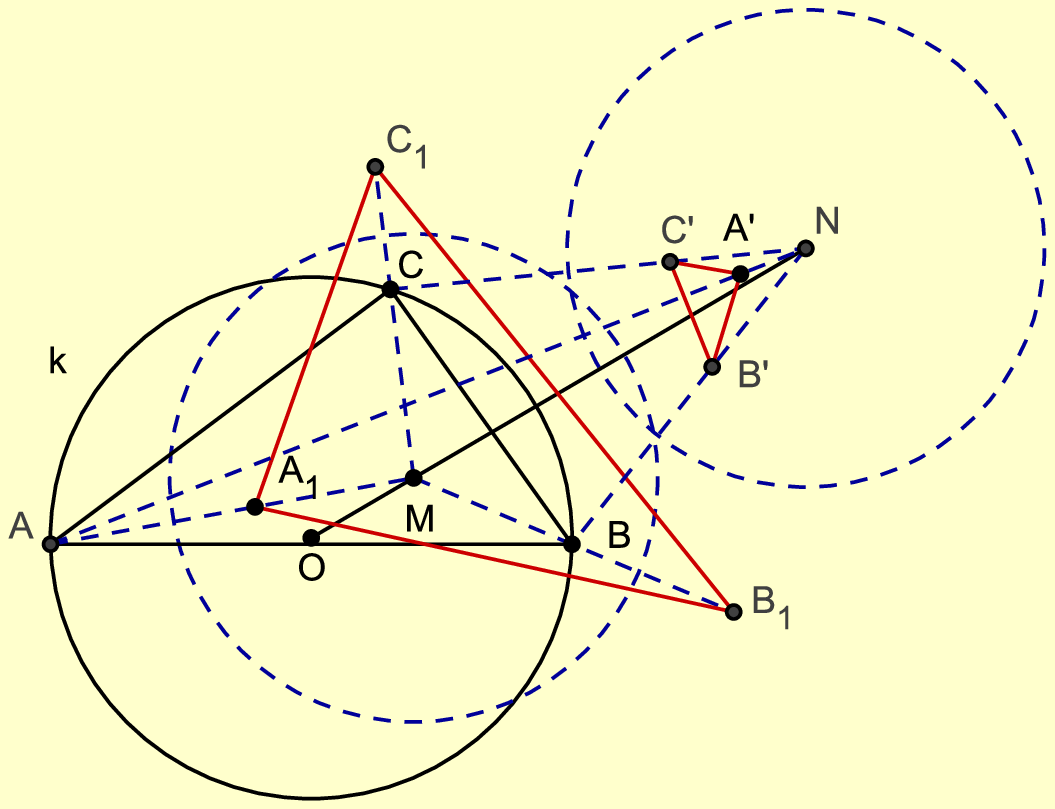}{Figure 1}
\end{figure}

The points $A_1,B_1,C_1$ are collinear if and only if $M \in k$.

Next we assume $M\notin k$ and denote the angles of $\triangle A_1B_1C_1$ by
$(\alpha_1, \beta_1, \gamma_1)$. Then (1) means that any circle inversion
$\varphi(M)$ with center $M$ determines one and the same triple of angles
$(\alpha_1, \beta_1, \gamma_1)$ of a triangle. Hence, there arises a map
$$f: \quad M \longrightarrow (\alpha_1, \beta_1, \gamma_1);
\qquad M \notin k, \quad (\alpha_1, \beta_1, \gamma_1) \, -
{\rm a \, triple \, of \, angles \, of \, a \, triangle}.$$

First we note that $f(O)=(\alpha, \beta, \gamma)$.

For points, different from $O$, we shall prove the following statement.

\begin{prop}\label{P:2.1}
Given a $\triangle ABC$ and the circum-circle $k$ of the triangle.
If $M$ and $N$ are two inverse points with respect to the circle $k$,
then
$$f(M) = f(N).$$
\end{prop}

\emph{Proof:}  Let $\varphi(M,r)$ and $\varphi'(N,r')$ be two circle inversions
with centers $M$ and $N$, respectively. Denote by $A_1B_1C_1$ the image
of the triple $ABC$ under the inversion $\varphi$ and by $A'B'C'$ the image
of the triple $ABC$ under the inversion $\varphi'$ (Figure 1).

If $MA=x$, $MB=y$, $MC=z$ and $NA=x', NB=y', NC=z'$, then we have
$$\begin{array}{l}
B_1C_1=\ds{\frac{r^2}{xyz}\,ax, \; C_1A_1=\frac{r^2}{xyz}\,by, \;
A_1B_1=\frac{r^2}{xyz}\,cz;}\\
[4mm]
B'C'=\ds{\frac{r'^2}{x'y'z'}\,ax', \; C'A'=\frac{r'^2}{x'y'z'}\,by',
\; A'B'=\frac{r'^2}{x'y'z'}\,cz'.} \end{array}\leqno(2)$$

On the other hand, since $M$ and $N$ are inverse points with respect
to the circle $k$, then $k$ is Apollonius circle with basic points $M,N$
and ratio
$\ds{\frac{d}{R}}$, where $OM=d$. Therefore
$$\frac{x}{x'}=\frac{y}{y'}=\frac{z}{z'}=\frac{d}{R}\, .\leqno(3)$$
Then (2) and (3) imply that
$$\frac{B_1C_1}{B'C'}=\frac{C_1A_1}{C'A'}=\frac{A_1B_1}{A'B'}
=\frac{r^2}{r'^2}\,\frac{x'y'z'}{xyz}\,
\frac{d}{R}.$$
Hence $\triangle A_1B_1C_1 \sim \triangle A'B'C'$, which means that
$f(M)=f(N)$. \qed
\begin{prop}\label{P:2.2}
Given the basic triangle $ABC$ with angles $(\alpha, \beta, \gamma)$
and an arbitrary triple of angles $(\alpha_1, \beta_1, \gamma_1)$ of a triangle.
Prove that:

(i) if $(\alpha_1, \beta_1, \gamma_1) \neq (\alpha, \beta, \gamma)$,
then there exist exactly two points $M$ and $N$ such that
$f(M)=f(N)=(\alpha_1, \beta_1, \gamma_1)$;

(ii) if $(\alpha_1, \beta_1, \gamma_1) = (\alpha, \beta, \gamma)$,
then only the center $O$ satisfies the condition $f(O)=(\alpha, \beta, \gamma)$.
\end{prop}

\emph{Proof:} (i) Let $A_1B_1C_1$ be a triangle with angles
$(\alpha_1, \beta_1, \gamma_1)$, respectively, and denote by
$B_1C_1=a_1, \, C_1A_1=b_1, \, A_1B_1=c_1$. If $M$ is a solution
to the equation $f(X)=(\alpha_1, \beta_1, \gamma_1)$, then according to (1)
we have
$$\frac{ax}{a_1}=\frac{by}{b_1}=\frac{cz}{c_1}$$
and consequently
$$\frac{MB}{MC}=\frac{y}{z}=\frac{b_1}{c_1}\,\frac{c}{b}, \quad
\frac{MC}{MA}=\frac{z}{x}=\frac{c_1}{a_1}\,\frac{a}{c}, \quad
\frac{MA}{MB}=\frac{x}{y}=\frac{a_1}{b_1}\,\frac{b}{a}.$$

Thus we obtained that any solution $M$ to the equation
$f(X)=(\alpha_1, \beta_1, \gamma_1)$
is a common point of the three Apollonius circles
$k_1\left(B, C; \; \lambda=\ds{\frac{b_1}{c_1}\,\frac{c}{b}}\right)$,
$k_2\left(C, A; \; \mu=\ds{\frac{c_1}{a_1} \,\frac{a}{c}}\right)$
and \\ $k_3\left(A, B; \; \nu=\ds{\frac{a_1}{b_1}\,\frac{b}{a}}\right)$
(Figure 2).

\begin{figure}[h]\center\epsfysize=5cm\epsffile{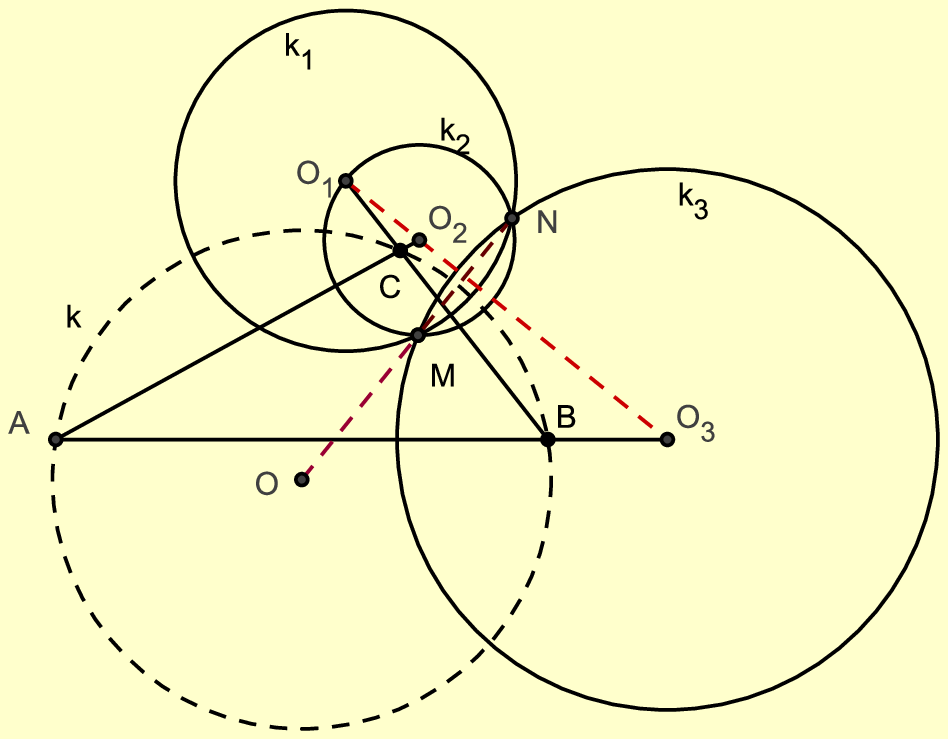}
{Figure 2}
\end{figure}

Since $(\alpha_1, \beta_1, \gamma_1) \neq (\alpha, \beta, \gamma)$
and $\lambda\mu\nu=1$, at least one of the numbers $\lambda, \mu, \nu$
is greater than $1$ and at least one of them is less than $1$.
For definiteness let $\lambda >1$ and $\mu<1$. Therefore the point
$C$ is interior for the circles $k_1$ and $k_2$.

Let $(O_1,r_1)$ and $(O_2,r_2)$ be the corresponding centers and radii
of $k_1$ and $k_2$, respectively. To prove that $k_1$ and $k_2$ intersect
into two points, it is sufficient to show that $|r_1-r_2|<O_1O_2$, which is
equivalent to
$$O_1O_2^2-(r_1-r_2)^2>0.$$

It is easy to find that
$$\overrightarrow{O_1O_2}=\frac{\mu^2}{1-\mu^2} \, \overrightarrow{BC}
-\frac{1}{\lambda^2-1} \, \overrightarrow{AC}, \quad
r_1=\frac{\lambda \, a}{\lambda^2-1}, \quad r_2=\frac{\mu \, b}{1-\mu^2}.$$
Then we get
$$O_1O_2^2-(r_1-r_2)^2=\frac{1}{(\lambda^2-1)(1-\mu^2)}
\,[(c\mu)^2-(b\lambda\mu-a)^2].$$
The equalities $\lambda=\ds{\frac{b_1}{c_1}\,\frac{c}{b}}$ and
$\mu=\ds{\frac{c_1}{a_1}\,\frac{a}{c}}$ imply that
$$\left|\frac{b\lambda}{c}-\frac{a}{cm}\right|=\frac{|b_1-a_1|}{c_1}<1.$$
Therefore $(c\mu)^2-(b\lambda\mu-a)^2>0$ and the circles $k_1$, $k_2$
intersect into two points, which we denote by $M$ and $N$.

The condition $\lambda \mu \nu = 1$ implies that the third circle $k_3$
also passes through $M$ and $N$.

(ii) If $\alpha_1=\alpha, \, \beta_1=\beta, \, \gamma_1=\gamma$, then
$\lambda = \mu = \nu =1$ and $k_1, \, k_2, \, k_3$ are
the perpendicular bisectors of the corresponding sides of $\triangle ABC$.
Hence, in this case $k_1, \, k_2$ and $k_3$ have one common point $O$ and
the only solution to the equation
$f(X)=(\alpha_1, \beta_1, \gamma_1)$ is the point $O$. \qed
\vskip 2mm
Now Proposition \ref{P:2.1} and Proposition \ref{P:2.2} imply the following

\begin{thm}\label{T:2.1}
Given the basic triangle $ABC$ with angles $(\alpha, \beta, \gamma)$ and
an arbitrary triple of angles $(\alpha_1, \beta_1, \gamma_1)$ of a triangle.
Then:

(i) there exists exactly one point $M$, interior for $k$, such that
$f(M)=(\alpha_1, \beta_1, \gamma_1)$;

(ii) if $(\alpha_1, \beta_1, \gamma_1) \neq (\alpha, \beta, \gamma)$,
then there exists exactly one point $N$, exterior for $k$, such that
$f(N)=(\alpha_1, \beta_1, \gamma_1)$. Furthermore $N$ and $M$ from (i)
are inverse points with respect to the circum-circle of $\triangle ABC$.
\end{thm}

\subsection{A realization of the map $f$}

Let $\triangle ABC$ be the basic triangle. For an arbitrary point $M$,
interior for the circum-circle $k(O,R)$ we denote by $A_1, B_1, C_1$
the orthogonal projections of $M$ on the lines $BC, CA, AB$, respectively.
Then $\triangle A_1B_1C_1$ is the \emph{pedal triangle} of $M$ with respect
to $\triangle ABC$. Denoting by $\alpha_1, \beta_1, \gamma_1$, respectively,
the angles of $\triangle A_1B_1C_1$, then we can consider the correspondence
$$M \quad \longrightarrow \quad (\alpha_1, \beta_1, \gamma_1).$$

Next we shall prove that this correspondence is a realization of the map $f$.

\begin{prop} \label{P:3.1}
Given the basic $\triangle ABC$ with circum-circle $k$ and a point $M$,
interior for the circle $k$. If $(\alpha_1, \beta_1, \gamma_1)$ are
the corresponding angles of the pedal triangle of $M$ with respect to
the basic triangle, then
$$f(M)=(\alpha_1, \beta_1, \gamma_1).$$
\end{prop}

\emph{Proof:} Since $AM=x$ is a diameter of the circum-circle of
$\triangle AC_1B_1$ (Figure 3), we find
$$B_1C_1=AM \sin \alpha= \frac{ax}{2R}.$$
Similarly we find $C_1A_1=\ds{\frac{by}{2R}}$ and $A_1B_1={\ds \frac{cz}{2R}}$.
Thus we obtain
$$\frac{ax}{B_1C_1}=\frac{by}{C_1A_1}=\frac{cz}{A_1B_1}=2R.$$
Now the assertion follows from (1) and the definition of $f$. \qed

\begin{figure}[h]\center\epsfysize=5cm\epsffile{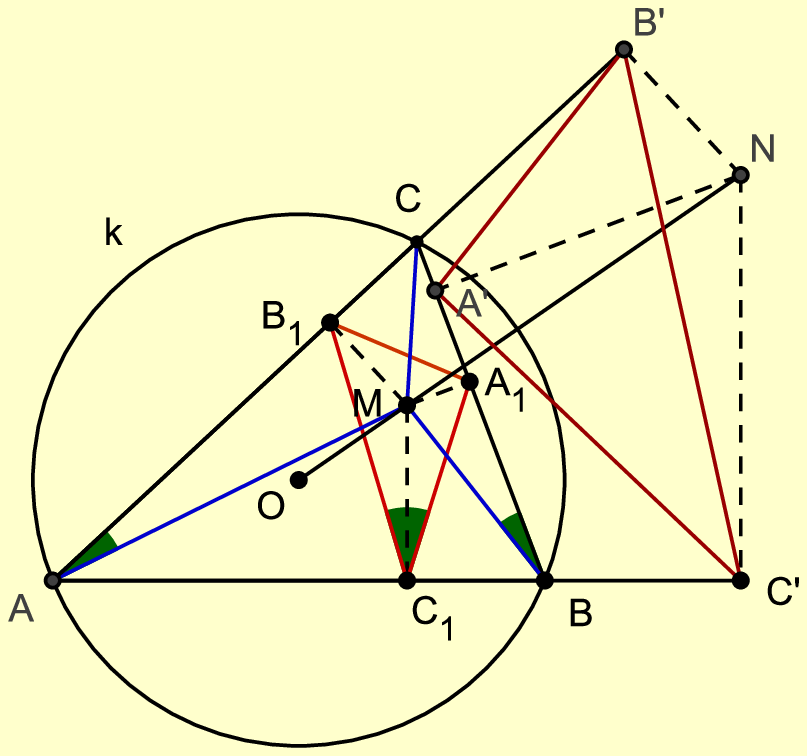}
{Figure 3}
\end{figure}

Taking into account how a circle inversion transforms the circles in the plane,
we have:

{\it Let $\varphi$ be a circle inversion with center $M$ and $\varphi(ABC)=A'B'C'$.
Then:
\begin{itemize}
\item $\triangle A'B'C'$ is positively oriented (has the orientation of
$\triangle ABC$) if and only if $M$ is inside the circle $k$;
\item $\triangle A'B'C'$ is negatively oriented if and only if $M$ is
outside the circle $k$;
\end{itemize}}

The above statement and Proposition 2.1 imply the following:

{\it Let $\triangle A_1B_1C_1$ be the pedal triangle of a point $M$.
Then:
\begin{itemize}
\item $\triangle A_1B_1C_1$ is positively oriented (has the orientation of
$\triangle ABC$) if and only if $M$ is inside the circle $k$;
\item $\triangle A_1B_1C_1$ is negatively oriented if and only if $M$ is
outside the circle $k$;
\end{itemize}}

A point $M$ lies on $k$, if and only if its pedal triangle is degenerate
(which is the Simson's theorem).

These considerations, Proposition \ref{P:3.1} and Theorem \ref{T:2.1} imply
the following

\begin{thm}\label{T:3.1}
Given the basic triangle $ABC$ with angles $(\alpha, \beta, \gamma)$ and
an arbitrary triple of angles $(\alpha_1, \beta_1, \gamma_1)$ of a triangle.
Then:

(i) there exists exactly one point $M$, inside the circum-circle $k$ of
$\triangle ABC$, such that $(\alpha_1, \beta_1, \gamma_1)$ are the angles
of its pedal triangle. In this case the pedal triangle of $M$ and the basic
triangle have the same orientation.

(ii) if $(\alpha_1, \beta_1, \gamma_1) \neq (\alpha, \beta, \gamma)$,
then there exists exactly one point $N$, outside the circum-circle $k$
of $\triangle ABC$, such that $(\alpha_1, \beta_1, \gamma_1)$ are the angles
of its pedal triangle. In this case the pedal triangle of $N$ and the basic
triangle have the opposite orientations. Furthermore, $N$ and $M$ from (i)
are inverse points with respect to the circum-circle of the basic triangle.
\end{thm}

\begin{rem}
It follows from Theorem \ref{T:3.1} that there does not exist a point $N$,
whose pedal triangle is similar to the basic triangle and negatively oriented.
\end{rem}

The angles of $\triangle ABC$ and $\triangle A_1B_1C_1$ determine the angles
$\angle BMC, \angle CMA$, and $\angle AMB$ in the following way (Figure 3):
$$\angle AMB=\gamma + \gamma_1, \quad \angle CMA=\beta + \beta_1, \quad
\angle BMC=\alpha + \alpha_1.\leqno(1)$$

We adopt the following convention: $\angle BMC=\pi$ if and only if
$M$ lies on the side $BC$; $\angle BMC>\pi$ if and only if $M$ is inside
the circle $k$, but $M$ and $A$ are from different sides of
the line $BC$. Then formulas (1) are valid for all points $M$ inside
the circle $k$. Then, for any point inside the circle $k$, we have
the following simple criterion:

\emph{Let $M$ be inside the circle $k$. Then the angles of
the pedal $\triangle A_1B_1C_1$ are
$$\alpha_1= \angle BMC - \alpha, \quad \beta_1 = \angle CMA - \beta, \quad
\gamma_1 = \angle AMB - \gamma.$$}

\section{Points, whose pedal triangles are similar to the given one}
Let $k(O,R)$ be the circum-circle of the basic $\triangle ABC$. In this section
we study the points, whose pedal triangles are similar to the given $\triangle ABC$.

\subsection{Points, whose pedal triangles are positively oriented}
Let $\Omega_1$ be the first Brocard point for the basic $\triangle ABC$
(Figure 4). According to the definition of $\Omega_1$ we have
$\angle \Omega_1AB=\angle \Omega_1BC=\angle \Omega_1CA=\omega$, $\omega$ being the
Brocard angle.

\begin{figure}[h]\center\epsfysize=5cm\epsffile{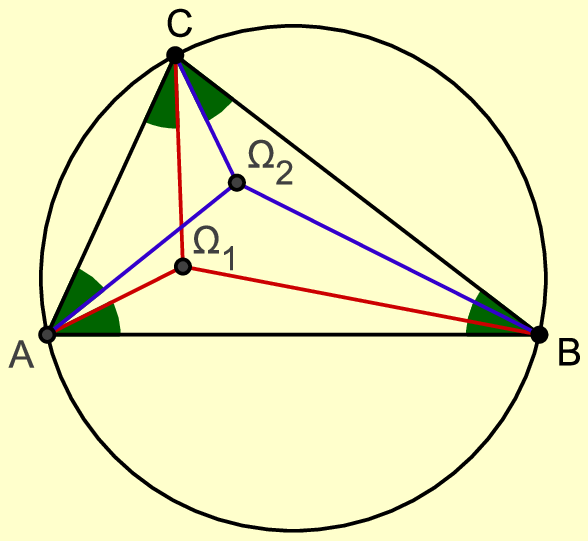}
{Figure 4}
\end{figure}

If $A_1B_1C_1$ is the pedal triangle of $\Omega_1$, then it follows
immediately that $\angle B_1A_1C_1=\alpha_1=\beta, \;
\angle C_1B_1A_1=\beta_1=\gamma, \; \angle A_1C_1B_1= \gamma_1=\beta$.
Hence, the pedal triangle of $\Omega_1$ is similar to $\triangle BCA$
and positively oriented.

The second Brocard point $\Omega_2$ is introduced by the conditions
$\angle \Omega_2AC=\angle \Omega_2CB=\angle \Omega_2BA=\omega$,
which implies that the angles of the pedal triangle of $\Omega_2$ are
$(\alpha_1,\beta_1,\gamma_1)=(\gamma, \alpha, \beta)$.
Therefore, the pedal triangle of $\Omega_2$ is similar
to $\triangle CAB$ and positively oriented.

Now, let $L_3$ be the point, whose pedal $\triangle A_1B_1C_1$ is positively
oriented and similar to $\triangle BAC$. This means that $L_3$ is inside
the circle $k$ and
$$\angle BL_3C=\angle CL_3A=\alpha+\beta, \qquad \angle AL_3B=2\gamma.$$
The last equality implies that $L_3$ is a point on the arc $AOB$ (Figure 5).

Further, we denote by $k'$ the circum-circle of $\triangle ABO$ ($k'$ is
the line $AB$ if $\angle ACB = 90^0$) and by
$EF$ the diameter of $k$, perpendicular to the side $AB$. If $K$ is
the common point of $CL_3$ and $EF$, then the condition
$\angle AL_3K=\angle BL_3K$ implies that $OK$ is a diameter of $k'$.

Let $\varphi (O, R)$ be the inversion with respect to the circle $k$.
If $D$ is the midpoint of the side $AB$, then $\varphi(D)=K$ and the quadruple
$DKEF$ is harmonic. Therefore, the circle $k$ is Apollonius circle
with basic points $D,K$ passing through the point $C$. Hence, $CE \,(CF)$ is the
interior (exterior) bisector of $\angle DCK$. Thus we obtained that $L_3$
is the common point of the arc $AOC$ and the ray of the symmedian to
the side $AB$. Since $OK$ is a diameter of the circle $k'$, then
$$OL_3 \perp CK.\leqno(2.1)$$

\begin{figure}[h]\center\epsfysize=8cm\epsffile{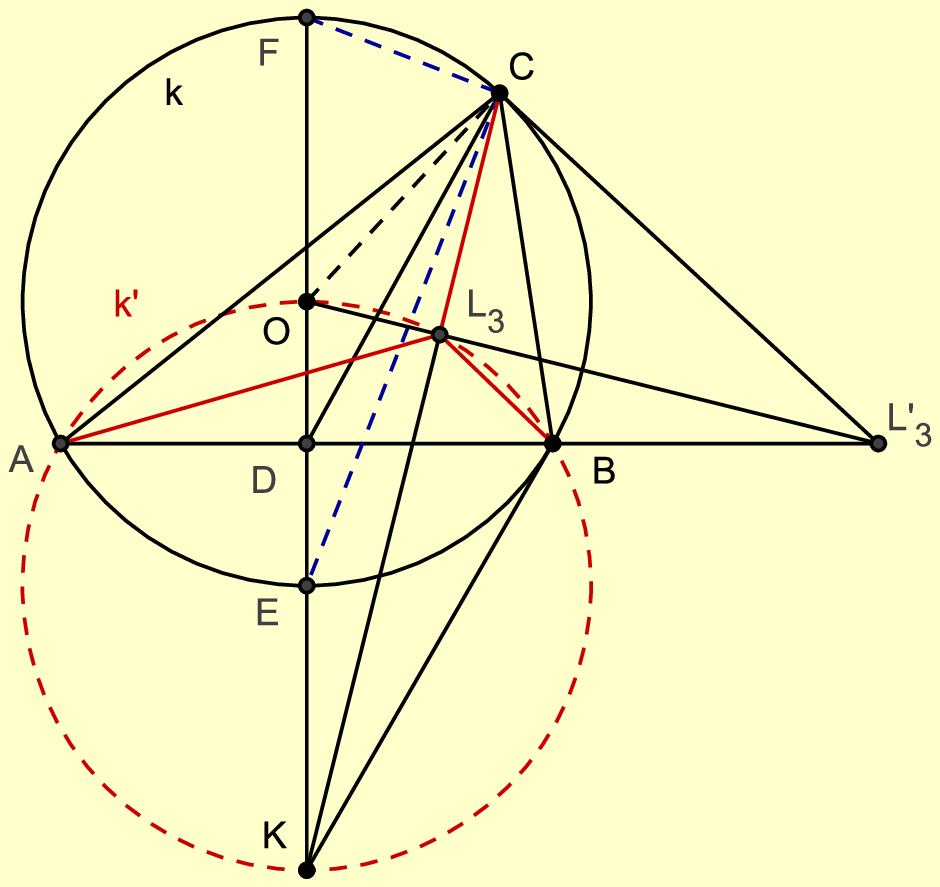}
{Figure 5}
\end{figure}

The orthogonal projections of the circumcenter $Î$ on the three symmedians of
$\triangle ABC$ are considered from many other aspects (see e.g. \cite{C}, Theorems
120-121).

So, the six points on the circle $k$, whose pedal triangles
are positively oriented and similar to the basic triangle, are as follows:

\begin{itemize}
\item
the center $O$ of the circum-circle: $\triangle A_1B_1C_1 \sim \triangle ABC$;
\item
the first Brocard point $\Omega_1$: $\triangle A_1B_1C_1 \sim \triangle BCA$;
\item
the second Brocard point $\Omega_2$: $\triangle A_1B_1C_1 \sim \triangle CAB$;
\item
the orthogonal projection $L_1$ of $O$ on the ray of the symmedian to the side $BC$:
$\triangle A_1B_1C_1 \sim \triangle ACB$
\item
the orthogonal projection $L_2$ of $O$ on the ray of the symmedian to the side $CA$:
$\triangle A_1B_1C_1 \sim \triangle CBA$;
\item
the orthogonal projection $L_3$ of $O$ on the ray of the symmedian to the side $AB$:
$\triangle A_1B_1C_1 \sim \triangle BAC$;
\end{itemize}

Let us denote by $L$ the Lemoine point for the basic triangle. It is wellknown fact
(e.g. \cite{C, H}) that the Brocard points $\Omega_1$ and $\Omega_2$ lie on the circle
$k_0$ with diameter $OL$ so that $\angle \Omega_1OL=\angle \Omega_2OL=\omega$, where
$\omega$ is the Brocard angle.

The circle $k_0$ is called \emph{the Brocard's circle}.

Now the following statement follows at once.

\begin{thm} The six points $O, \Omega_1, \Omega_2, L_1, L_2, L_3$, whose pedal triangles
are positively oriented and similar to the basic triangle, lie on the Brocard's circle.
\end{thm}

Using the properties of the symmedians of a triangle, we find the
following representations for the points $L_1, L_2, L_3$ in barycentric
coordinates:
$$\overrightarrow{OL}_1=\frac{(b^2+c^2-a^2)\,\overrightarrow{OA}
+b^2\,\overrightarrow{OB}+c^2\,\overrightarrow{OC}}{2(b^2+c^2)-a^2}\,,
\; \overrightarrow{OL}_2=\frac{a^2\,\overrightarrow{OA}
+(c^2+a^2-b^2)\,\overrightarrow{OB}+c^2\,\overrightarrow{OC}}{2(c^2+a^2)-b^2}\,,$$
$$\overrightarrow{OL}_3=\frac{a^2\,\overrightarrow{OA}
+b^2\,\overrightarrow{OB}+(a^2+b^2-c^2)\,\overrightarrow{OC}}{2(a^2+b^2)-c^2}\,.
$$

\subsection{Points, whose pedal triangles are negatively oriented}
In this subsection we study the points exterior for the circle $k$, whose
pedal triangles are similar to the basic triangle. It follows from
Theorem 1.5 that these points are the inverse points of
$\Omega_1, \Omega_2, L_1, L_2, L_3$ with respect to $\varphi$ and we denote them
by $\Omega'_1, \Omega'_2, L'_1, L'_2, L'_3$, respectively.

\begin{figure}[h]\center\epsfysize=8cm\epsffile{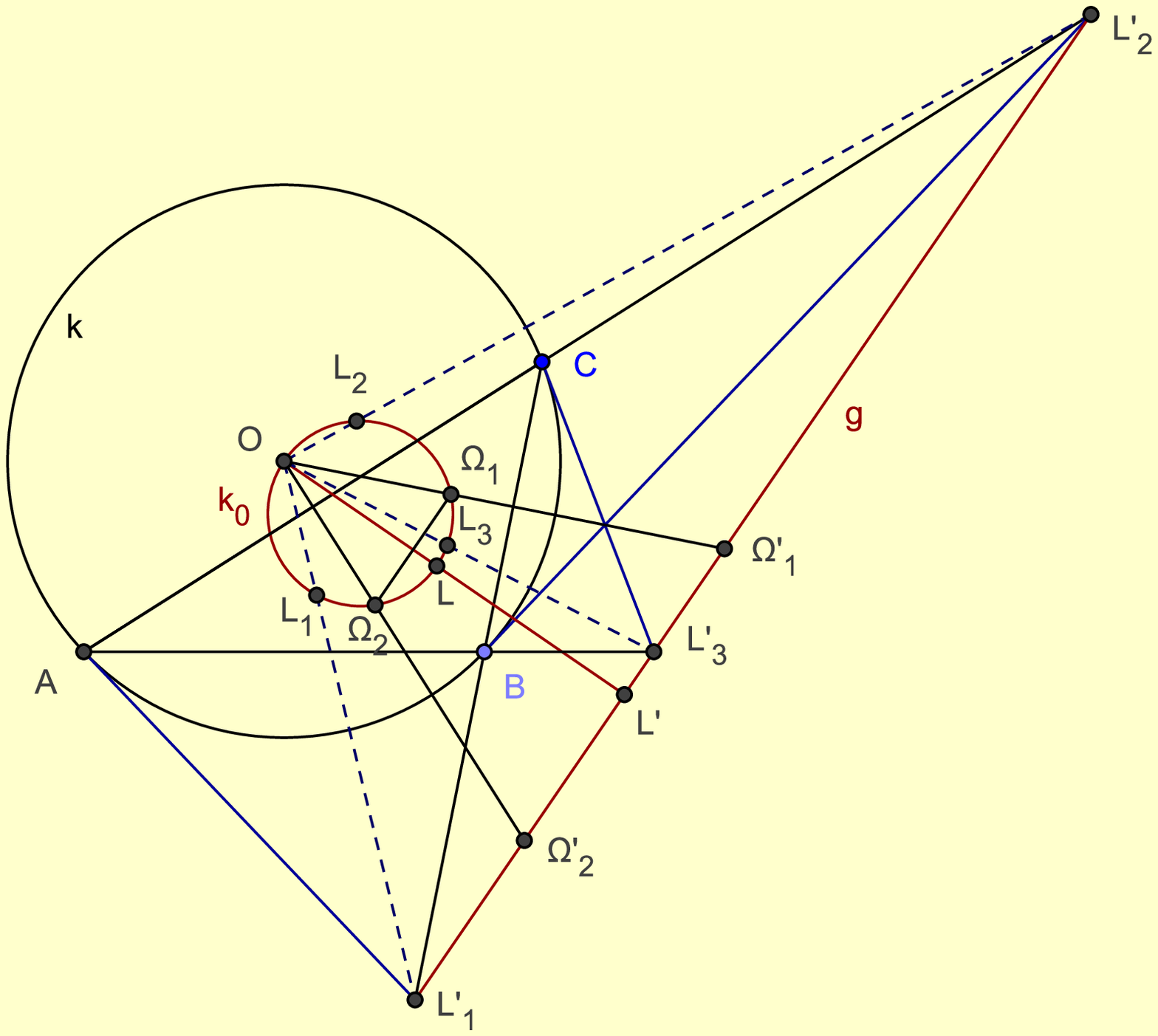}
{Figure 6}
\end{figure}

Applying Theorem 2.2 we obtain the following statement.

\begin{thm} The points $\Omega'_1, \Omega'_2, L'_1, L'_2, L'_3$,
whose pedal triangles are negatively oriented and similar to the basic
triangle, lie on a single straight line.
\end{thm}

We denote this straight line $g=\varphi(k_0)$.

The position of the five points from Theorem 2.2 is given in Figure 6 under
the assumption $a < c < b$.

First we study the point $L'_3$ (Figure 6). From one hand, $L'_3$ lies on
the line $OL_3$. On the other hand, $L'_3$ lies on the line $AB=\varphi(k')$.
Hence $L'_3=AB \cap OL_3$.

We note that if $CA=CB$ (isosceles $\triangle ABC$), then $g \Vert AB$
and $L'_3$ is the point at infinity of the line $AB$.

Taking into account that $L'_3=\varphi(L_3)$ and (2.1), we conclude
that $CL'_3$ is tangent to $k$ at $C$. Therefore $CL'_3$ is the exterior
simmedian through the vertex $C$ and $L'_3A:L'_3B=b^2:a^2$. The last
equality shows that the point $L'_3$ is the center of the Apollonius
circle $k_3$ with basic points $A,B$, which contains $C$. This circle becomes
the perpendicular bisector of the side $AB$ when $CA=CB$.

Similarly we have:
$L'_1$ is the center of the Apollonius circle $k_1$ with
basic points $B,C$, containing $A$; $L'_2$ is the center of the Apollonius
circle $k_2$ with basic points $C,A$, containing $B$. We call these Apollonius
circles \emph{the basic Apollonius circles} of $\triangle ABC$.

Further, we obtain the following statement.
\begin{cor} For any triangle the straight line, passing through the center
of the circum-circle and the Lemoine point, is perpendicular to the axis,
containing the centers of the three basic Apollonius circles of the triangle.
\end{cor}

The above statement also shows that $\Omega_1\Omega_2$ is parallel to
the axis $g$.

The points $\Omega'_1$ and $\Omega'_2$ are $O\Omega_1 \cap g$ and
$O\Omega_2 \cap g$, respectively, and $L'$ is the midpoint of the segment
$\Omega'_1\Omega'_2$. Then we can calculate the barycentric coordinates
of the point $\Omega_1'$:

$$O\Omega'_1=\frac{a^2(a^2-b^2)\,\overrightarrow{OA}+b^2(b^2-c^2)\,\overrightarrow{OB}
+c^2(c^2-a^2)\,\overrightarrow{OC}}{a^4+b^4+c^4-a^2b^2-b^2c^2-c^2a^2}\,.$$
Replacing $(a,b,c)$ with $(b,c,a)$ in the last formula, we obtain
$$O\Omega'_2=\frac{a^2(a^2-c^2)\,\overrightarrow{OA}+b^2(b^2-a^2)\,\overrightarrow{OB}
+c^2(c^2-b^2)\,\overrightarrow{OC}}{a^4+b^4+c^4-a^2b^2-b^2c^2-c^2a^2}\,.$$

Since $L'_3$ is the midpoint of the segment $\Omega'_1\Omega'_2$, then
$$OL'=\frac{a^2(2a^2-b^2-c^2)\,\overrightarrow{OA}+b^2(2b^2-c^2-a^2)\,\overrightarrow{OB}
+(2c^2-a^2-b^2)\,\overrightarrow{OC}}{2(a^4+b^4+c^4-a^2b^2-b^2c^2-c^2a^2)}\,.$$
\vskip 2mm
\begin{rem}
Taking into account the characterization:

{\it The first (second) Brocard point $\Omega_1$ ($\Omega_2$) is the only point,
whose pedal triangle is positively oriented and similar to $\triangle BCA$
($\triangle CAB$),}

\noindent
we can consider the points $\Omega'_1$ and $\Omega'_2$ as the
{\it exterior Brocard points} of the given triangle.
\end{rem}

From the above formulas we obtain immediately:
\begin{cor}
The following statements are equivalent:

(i) \; $b=c$;

(ii) \; $\Omega_1\equiv L_2$;

(iii) \; $\Omega_2\equiv L_3.$
\end{cor}
Further we obtain the following new characterization for the triangles, whose sides are
proportional to their medians.

\begin{prop}
Given a $\triangle ABC$ with circum-circle $k(O)$. The following conditions are equivalent:

(i) \; $a^2+b^2=2c^2$;

(ii) \; the Lemoine point $L$ lies on the circum-circle of $\triangle ABO$;

(iii) \; the inverse image of the Lemoine point with respect to $k(O)$ lies on
the side-line $AB$.
\end{prop}

\end{document}